# A CONNECTION BETWEEN THE GHIRLANDA–GUERRA IDENTITIES AND ULTRAMETRICITY[1]

By Dmitry Panchenko

*Texas A&M University*

We consider a symmetric positive definite weakly exchangeable infinite random matrix and show that, under the technical condition that its elements take a finite number of values, the Ghirlanda–Guerra identities imply ultrametricity.

**1. Introduction and main result.** Let us consider an infinite random matrix $R = (R_{l,l'})_{l,l' \geq 1}$ which is symmetric, nonnegative definite [in the sense that $(R_{l,l'})_{1 \leq l,l' \leq n}$ is nonnegative definite for any $n \geq 1$] and weakly exchangeable, which means that for any $n \geq 1$ and any permutation $\rho$ of $\{1, \ldots, n\}$, the matrix $(R_{\rho(l),\rho(l')})_{1 \leq l,l' \leq n}$ has the same distribution as $(R_{l,l'})_{1 \leq l,l' \leq n}$. Following [6], we will call the matrix with such properties a Gram–de Finetti matrix. We assume that diagonal elements $R_{l,l} = 1$ and nondiagonal elements take only a finite number of values,

$$(1.1) \qquad \mathbb{P}(R_{1,2} = q_l) = m_{l+1} - m_l$$

for $1 \leq l \leq k$ and for some $-1 \leq q_1 < q_2 < \cdots < q_k \leq 1$ and $0 = m_1 < \cdots < m_k < m_{k+1} = 1$. We say that the matrix $R$ satisfies the *Ghirlanda–Guerra identities* [7] if, for any $n \geq 2$, any bounded measurable functions $f : \mathbb{R}^{n(n-1)/2} \to \mathbb{R}$ and $\psi : \mathbb{R} \to \mathbb{R}$,

$$(1.2) \qquad \mathbb{E} f_n \psi(R_{1,n+1}) = \frac{1}{n} \mathbb{E} f_n \mathbb{E} \psi(R_{1,2}) + \frac{1}{n} \sum_{l=2}^{n} \mathbb{E} f_n \psi(R_{1,l}),$$

where $f_n = f((R_{l,l'})_{1 \leq l < l' \leq n})$. In other words, conditionally on $(R_{l,l'})_{1 \leq l < l' \leq n}$, the law of $R_{1,n+1}$ is given by the mixture $n^{-1} \mathcal{L}(R_{1,2}) + n^{-1} \sum_{l=2}^{n} \delta_{R_{1,l}}$. By the positivity principle proven by Talagrand (Theorem 6.6.2 in [14], see also

Received March 2009.
[1]Supported in part by NSF grant DMS-0904565.
*AMS 2000 subject classifications.* 60K35, 82B44.
*Key words and phrases.* Sherrington–Kirkpatrick model, Parisi ultrametricity conjecture.

[9]), the Ghirlanda–Guerra identities imply that $R_{1,2} \geq 0$ with probability 1 and, therefore, from now on, we can assume that $q_1 \geq 0$. However, this a priori assumption is not necessary since it will also be clear from the proof that $q_1$ must be nonnegative. The main result of the paper is the following.

THEOREM 1.    *Under assumptions (1.1) and (1.2), the matrix $R$ is ultrametric, that is,*

$$(1.3) \qquad \mathbb{P}(R_{2,3} \geq \min(R_{1,2}, R_{1,3})) = 1.$$

Another way to express the event in (1.3) is to say that

$$(1.4) \qquad R_{1,2} \geq q_l, R_{1,3} \geq q_l \implies R_{2,3} \geq q_l \qquad \text{for all } 1 \leq l \leq k.$$

This question of ultrametricity originates in the setting of the Sherrington–Kirkpatrick model [13], where $R$ corresponds to the matrix of the overlaps, or scalar products, of i.i.d. replicas from a random Gibbs measure. The ultrametricity property (1.3) was famously predicted by Parisi in [8] as a part of complete description of the expected behavior of the model and it still remains an open mathematical problem. On the other hand, the Ghirlanda–Guerra identities, which are implicitly contained in the Parisi theory, were proven rigorously in [7] in some approximate sense; namely, one can slightly perturb the parameters of the model such that, on average over the perturbation (or for some specific choice of perturbed parameters, see [15]), the Ghirlanda–Guerra identities hold in the thermodynamic limit. In this paper, we consider an asymptotic distribution of the overlap matrix for which the Ghirlanda–Guerra identities hold precisely, as in (1.2), and show that they automatically imply ultrametricity, under a technical condition (1.1). In some sense, the main idea of the paper is nothing but a reversal of the proof of the Ghirlanda–Guerra identities, which arise from the information provided by the "stochastic stability" of the system with regard to small perturbations of the parameters. Our main technical contribution, the invariance principle of Theorem 4, is a very specific form of the stochastic stability of the system implied by the Ghirlanda–Guerra identities. The Parisi theory for the Sherrington–Kirkpatrick model also predicts that the distribution of the overlap $R_{1,2}$ has a nontrivial continuous component, so the assumption (1.1) is rather restrictive. Nevertheless, Theorem 1 provides some hope that the Ghirlanda–Guerra identities imply ultrametricity in the general case as well.

This work was motivated by a paper of Arguin and Aizenman [3] and, in particular, by a beautiful application of the Dovbysh–Sudakov representation theorem [6] for Gram–de Finetti matrices that will play the same crucial role here. In [3], the authors prove ultrametricity in a slightly different setting as a consequence of what they call the "robust quasi-stationarity property." The main idea in [3] utilizes the robust quasi-stationarity in order to prove "quasi-stationarity under free evolution" at each step of the inductive argument, which, in turn, implies weak exchangeability and, via the application of the Dovbysh–Sudakov representation, induces clustering of the type (1.4). Our proof is based on exactly the same idea. The difference now is that quasi-stationarity under free evolution—the invariance or stochastic stability property of Theorem 4 below—will be a consequence of the Ghirlanda–Guerra identities and, of course, the induction in the proof of Theorem 1 will be different since it is also based on (1.2). In addition, we give a new proof in Theorem 3 below that the invariance implies exchangeability, which is based on the explicit control of the mixing induced by the random permutation in the invariance principle.

Simultaneously with the present work, Talagrand developed a different approach to Theorem 1 in [15] based on the Ghirlanda–Guerra identities and a form of invariance. Theorem 4 below shows that sufficient invariance is already contained in the Ghirlanda–Guerra identities and one can now find a new more direct proof of Theorem 4 in [15]. In addition, [15] clarifies the physicists' idea of decomposing the system into pure states and explains how both the Ghirlanda–Guerra identities and invariance arise in the Sherrington–Kirkpatrick model.

The rest of the paper is organized as follows. In the next section, we start with the Dovbysh–Sudakov representation result which shows that any Gram–de Finetti matrix $R$ can be generated by i.i.d. replicas from some random Gibbs measure on a separable Hilbert space, which is called the *directing measure* of $R$, in almost exactly the same way as the overlap matrix is generated by the Gibbs measure in the Sherrington–Kirkpatrick model. We first study some basic properties of the directing measure which follow from the Ghirlanda–Guerra identities; namely, that it always concentrates on a nonrandom sphere of the Hilbert space and is either continuous or discrete with probability 1. In particular, it is discrete when (1.1) holds. In Section 3, we formulate the invariance and exchangeability properties of the configuration of the directing measure and use them to prove Theorem 1 by induction on $k$. Finally, in Section 4, we show how the Ghirlanda–Guerra identities imply the invariance of the directing measure and how the invariance implies weak exchangeability.

**2. Basic consequences of GGI and exchangeability.** Since all of the properties of the matrix $R = (R_{l,l'})_{l,l' \geq 1}$ considered above—symmetry, positive definiteness, weak exchangeability, satisfying the Ghirlanda–Guerra identities (GGI)—were expressed in terms of its finite-dimensional distributions, we can think of $R$ as a random element in the product space $M = \prod_{1 \leq l,l'}[-1,1]$ with the pointwise convergence topology and the Borel $\sigma$-algebra $\mathcal{M}$. Let $\mathcal{P}$ denote the set of all probability measures on $\mathcal{M}$. Suppose that $\mathbb{P} \in \mathcal{P}$ is such that for all $A \in \mathcal{M}$,

$$(2.1) \qquad \mathbb{P}(A) = \int_\Omega \mathbb{Q}(\omega, A) \, d\Pr(\omega),$$

where $\mathbb{Q}: \Omega \times \mathcal{M} \to [0,1]$ is a probability kernel from some probability space $(\Omega, \mathcal{F}, \Pr)$ to $M$ such that: (a) $\mathbb{Q}(\omega, \cdot) \in \mathcal{P}$ for all $\omega \in \Omega$; (b) $\mathbb{Q}(\cdot, A)$ is measurable on $\mathcal{F}$ for all $A \in \mathcal{M}$. In this case, we will say that $\mathbb{P}$ is a *mixture* of laws $\mathbb{Q}(\omega, \cdot)$. Under (2.1), we can write the expectation of any measurable $\mathbb{P}$-integrable function $\phi: M \to \mathbb{R}$ as

$$(2.2) \qquad \mathbb{E}\phi(R) = \int_\Omega \left( \int_M \phi(R) \mathbb{Q}(\omega, dR) \right) d\Pr(\omega).$$

We will say that a law $\mathbb{Q} \in \mathcal{P}$ of a Gram–de Finetti matrix is *generated by an i.i.d. sample* if there exists a probability measure $\eta$ on $H \times [0, \infty)$, where $H$ is a separable Hilbert space, such that $\mathbb{Q}$ is the law of

$$(2.3) \qquad (x_l \cdot x_{l'} + a_l \delta_{l,l'})_{l,l' \geq 1},$$

where $(x_l, a_l)$ is an i.i.d. sequence from $\eta$ and $x \cdot y$ denotes the scalar product on $H$. The analysis of the distribution of $R$ will utilize the following representation result for Gram–de Finetti matrices due to Dovbysh and Sudakov [6].

PROPOSITION 1.  *A law $\mathbb{P} \in \mathcal{P}$ of any Gram–de Finetti matrix is a mixture (2.1) of laws in $\mathcal{P}$ such that for all $\omega \in \Omega$, $\mathbb{Q}(\omega, \cdot)$ is generated by an i.i.d. sample.*

We will denote by $\eta_\omega$ a probability measure on $H \times [0, \infty)$ corresponding to $\mathbb{Q}(\omega, \cdot)$ and let $\mu_\omega$ be the marginal of $\eta_\omega$ on $H$. Following the terminology of Aldous [2], we will call $\mu_\omega$ the *directing measure* of the matrix $(R_{l,l'})$. The main result of this section shows that the Ghirlanda–Guerra identities imply the following basic geometric properties of the directing measure.

THEOREM 2.  *Let $R$ be a Gram–de Finetti matrix on $M$ that satisfies the Ghirlanda–Guerra identities (1.2) and let $F$ be the law of $R_{1,2}$. If $q^*$ is the largest point of the support of $F$, then, for $\Pr$-almost all $\omega$: (a) $\mu_\omega(\|x\|^2 = q^*) = 1$; (b) $\mu_\omega$ is continuous if $F(\{q^*\}) = 0$; (c) $\mu_\omega$ is discrete if $F(\{q^*\}) > 0$, in which case the sequence of weights of $\mu_\omega$ has the Poisson–Dirichlet distribution $PD(1 - F(\{q^*\}))$.*

We recall that, given $s \in (0,1)$, if $(u_l)_{l \geq 1}$ is the decreasing enumeration of a Poisson point process on $(0, \infty)$ with intensity measure $x^{-1-s} \, dx$ and $w_l = u_l / \sum_j u_j$, then the distribution of the sequence $(w_l)$ is called the *Poisson–Dirichlet distribution* $PD(s)$. If $s = 0$, then we define $PD(0)$ to be the trivial distribution with $w_1 = 1$. The proof will be based on the following consequence of the Ghirlanda–Guerra identities.

LEMMA 1. *Consider a measurable set $A \subset [-1, 1]$. For $\Pr$-almost all $\omega$:*

(a) *if $F(A) > 0$, then for $\mu_\omega$-almost all $x_1$, $\mu_\omega(x_2 : x_1 \cdot x_2 \in A) > 0$;*
(b) *if $F(A) = 0$, then for $\mu_\omega$-almost all $x_1$, $\mu_\omega(x_2 : x_1 \cdot x_2 \in A) = 0$.*

PROOF. (a) Let $m = F(A^c) < 1$. Then, by the Ghirlanda–Guerra identities (1.2),

$$\mathbb{P}(\forall 2 \leq l \leq n+1, R_{1,l} \in A^c) = \mathbb{E} I(\forall 2 \leq l \leq n, R_{1,l} \in A^c) I(R_{1,n+1} \in A^c)$$

$$= \frac{n-1+m}{n} \mathbb{P}(\forall 2 \leq l \leq n, R_{1,l} \in A^c)$$

(by induction on $n$) $= \dfrac{(n-1+m)\cdots(1+m)m}{n!}$

$$= \frac{m(1+m)}{n}\left(1 + \frac{m}{2}\right)\cdots\left(1 + \frac{m}{n-1}\right)$$

(using $1 + x \leq e^x$) $\leq \dfrac{m(1+m)}{n} e^{m \log n} = \dfrac{m(1+m)}{n^{1-m}}.$

Since $m < 1$, letting $n \to +\infty$ gives $\mathbb{P}(\forall 2 \leq l, R_{1,l} \in A^c) = 0$. By Proposition 1, $\mathbb{P}$ is a mixture (2.1) of measures generated by an i.i.d. sample from $\eta_\omega$ and we can write

$$\mathbb{P}(\forall 2 \leq l, R_{1,l} \in A^c) = \int_\Omega \mu_\omega^{\otimes(l \geq 1)}(\forall 2 \leq l, x_1 \cdot x_l \in A^c) \, d\Pr(\omega) = 0,$$

which implies that for $\Pr$-almost all $\omega$, by Fubini's theorem,

$$\int \mu_\omega^{\otimes(l \geq 2)}((x_l)_{l \geq 2} : \forall l \geq 2, x_1 \cdot x_2 \in A^c) \, d\mu_\omega(x_1) = 0.$$

This implies that for $\mu_\omega$-almost all $x_1$,

$$\mu_\omega^{\otimes(l \geq 2)}((x_l)_{l \geq 2} : \forall l \geq 2, x_1 \cdot x_2 \in A^c) = \mu_\omega(x_2 : x_1 \cdot x_2 \in A^c)^\infty = 0,$$

which means that $\mu_\omega(x_2 : x_1 \cdot x_2 \in A^c) < 1$ and, therefore, $\mu_\omega(x_2 : x_1 \cdot x_2 \in A) > 0$. To prove part (b), it is enough to express $\mathbb{P}(R_{1,2} \in A) = 0$ using Proposition 1. □

PROOF OF THEOREM 2. (a) Since $F([-1, q^*]) = 1$ and $F([q^* - n^{-1}, q^*]) > 0$ for all $n \geq 1$, Lemma 1 implies that for $\Pr$-almost all $\omega$, for $\mu_\omega$-almost all $x_1$,

(2.4) $\quad \mu_\omega(x_2 : x_1 \cdot x_2 \leq q^*) = 1, \qquad \mu_\omega(x_2 : x_1 \cdot x_2 \geq q^* - n^{-1}) > 0$

$$\forall n \geq 1.$$

Let us fix any such $\omega$. First, the equality in (2.4) implies that $\mu_\omega(\|x\|^2 \leq q^*) = 1$. Otherwise, there exists $h \in H$ with $\|h\|^2 > q^*$ such that $\mu_\omega(B_\varepsilon(h)) > 0$ for any $\varepsilon > 0$, where $B_\varepsilon(h)$ is the ball of radius $\varepsilon > 0$ centered at $h$. Taking $\varepsilon > 0$ small enough so that $x_1 \cdot x_2 > q^*$ for all $x_1, x_2 \in B_\varepsilon(h)$ contradicts the first equality in (2.4). Next, let us show that $\mu_\omega(\|x\|^2 < q^*) = 0$. Otherwise, there again exists an open ball $B_\varepsilon(h) \subset \{x: \|x\|^2 < q^*\}$ such that $\mu_\omega(B_\varepsilon(h)) > 0$. For some $\delta > 0$, $\|x\|^2 < q^* - \delta$ for all $x \in B_\varepsilon(h)$ and, therefore, for all $x_1 \in B_\varepsilon(h)$ and $x_2 \in \{\|x\|^2 \leq q^*\}$, we have $x_1 \cdot x_2 < \sqrt{q^*(q^* - \delta)} \leq q^* - n^{-1}$ for large enough $n \geq 1$. Since we have already proven that $\mu_\omega(\|x\|^2 \leq q^*) = 1$, this contradicts the second inequality in (2.4). We have proven that $\mu_\omega(\|x\|^2 = q^*) = 1$.

(b) If $F(\{q^*\}) = 0$, then, by Lemma 1 for Pr-almost all $\omega$, for $\mu_\omega$-almost all $x_1$, we have $\mu_\omega(x_2 : x_1 \cdot x_2 = q^*) = 0$. Therefore, for Pr-almost all $\omega$ for which also $\mu_\omega(\|x\|^2 = q^*) = 1$, $\mu_\omega$ must be continuous.

(c) If $F(\{q^*\}) = 1$, then, by Lemma 1 for Pr-almost all $\omega$, for $\mu_\omega$-almost all $x_1$, we have $\mu_\omega(x_2 : x_1 \cdot x_2 = q^*) = 1$. By part (a), $\mu_\omega(\|x\|^2 = q^*) = 1$ and, therefore, $\mu_\omega$ must be concentrated on one point. If $F(\{q^*\}) \in (0,1)$, then for Pr-almost all $\omega$, for $\mu_\omega$-almost all $x_1$, we have $\mu_\omega(x_2 : x_1 \cdot x_2 = q^*) > 0$ and since $\mu_\omega(\|x\|^2 = q^*) = 1$, we get that for $\mu_\omega$-almost all $x_1$, $\mu_\omega(\{x_1\}) > 0$. This proves that $\mu_\omega$ is discrete. Let $(w_l)$ be the sequence of weights of $\mu_\omega$ arranged in decreasing order [we keep the dependence of $(w_l)$ on $\omega$ implicit]. The fact that $(w_l)$ has $PD(1 - F(\{q^*\}))$ distribution will follow from the analog of the Ghirlanda–Guerra identities for the Poisson–Dirichlet distribution proven by Talagrand in Chapter 1 of [14]. Let us explain how this result applies in our setting. Let us take any $m \geq 1$ and $n_1, \ldots, n_m \geq 1$ and for $n = n_1 + \cdots + n_m$, consider a function $f$ on $M$ which is the indicator of the set

$$\{R_{l,l'} = q^* : n_1 + \cdots + n_p + 1 \leq l \neq l' \leq n_1 + \cdots + n_p + n_{p+1}, 0 \leq p \leq m - 1\}.$$

Let us express the expectation $\mathbb{E}f$ using (2.2) and write the inside integral in terms of the weights $(w_l)$ of the directing measure $\mu_\omega$. Since, by part (a), $\mu_\omega(\|x\|^2 = q^*) = 1$ and $R_{l,l'} = x_l \cdot x_{l'} = q^*$ only if $x_1 = x_2$ when $\|x_1\|^2 = \|x_2\|^2 = q^*$, by Fubini's theorem,

$$(2.5) \qquad \int f(R)\mathbb{Q}(\omega, dR) = \sum_{l \geq 1} w_l^{n_1} \sum_{l \geq 1} w_l^{n_2} \cdots \sum_{l \geq 1} w_l^{n_m}.$$

Similarly, if we take $\psi(x) = I(x = q^*)$, then

$$(2.6) \qquad \int f\psi(R_{1,n+1})\mathbb{Q}(\omega, dR) = \sum_{l \geq 1} w_l^{n_1+1} \sum_{l \geq 1} w_l^{n_2} \cdots \sum_{l \geq 1} w_l^{n_m},$$

and we have $f\psi(R_{1,j}) = f$ for $2 \leq j \leq n_1$ and

$$(2.7) \quad \int f\psi(R_{1,j})\mathbb{Q}(\omega, dR) = \sum_{l \geq 1} w_l^{n_1+n_p} \cdots \sum_{l \geq 1} w_l^{n_{p-1}} \sum_{l \geq 1} w_l^{n_p+1} \cdots \sum_{l \geq 1} w_l^{n_m}$$

for $n_1 + 1 \leq j \leq n$ when $n_1 + \cdots + n_p + 1 \leq j \leq n_1 + \cdots + n_p + n_{p+1}$. If we let

$$S(n_1, \ldots, n_m) = \mathbb{E} \sum_{l \geq 1} w_l^{n_1} \sum_{l \geq 1} w_l^{n_2} \cdots \sum_{l \geq 1} w_l^{n_m}$$

and $s = 1 - F(\{q^*\})$, then plugging (2.5), (2.6) and (2.7) into (1.2) implies that

$$S(n_1 + 1, n_2, \ldots, n_m) = \frac{n_1 - s}{n} S(n_1, \ldots, n_m)$$
$$+ \sum_{p=2}^{m} \frac{n_p}{n} S(n_1 + n_p, \ldots, n_{p-1}, n_{p+1}, \ldots, n_m).$$

This coincides with equation (1.52) in [14]. It is explained there that this equation can be used recursively to compute $S(n_1, \ldots, n_m)$ in terms of $s \in (0, 1)$ only and that the set of numbers $S(n_1, \ldots, n_m)$ uniquely determines the distribution of the weights $(w_l)_{l \geq 1}$ as the Poisson–Dirichlet distribution $PD(s)$. □

**3. Ultrametricity in the discrete case.** In the case when (1.1) holds, Theorem 2 implies that for Pr-almost all $\omega$, the directing measure $\mu_\omega$ is discrete and concentrated on the sphere of radius $\sqrt{q_k}$, that is,

$$(3.1) \qquad \mu_\omega = \sum_{l \geq 1} w_l \delta_{\xi(l)}$$

for some distinct sequence $\xi(l) \in H$ with $\|\xi(l)\|^2 = q_k$ and $w_1 \geq w_2 \geq \cdots > 0$. In particular, by (2.3), this implies that a nondiagonal element $R_{l,l'} = x_l \cdot x_{l'} = q_k$ if and only if $x_l = x_{l'}$ which proves the *ultrametricity at the last level* $k$,

$$(3.2) \qquad R_{1,2} \geq q_k, \qquad R_{1,3} \geq q_k \implies R_{2,3} \geq q_k.$$

This is already a serious achievement, based, in addition to (1.1), only on the weak exchangeability of the matrix $R$ via the application of the Dovbysh–Sudakov representation. Bringing this idea to light was one of the main contributions of [3]. The description of the directing measure $\mu_\omega$ provided by (3.1) corresponds to the *pure states* picture in physics, where each point $\xi(l)$ can be thought of as a "pure state" of the asymptotic Gibbs measure. One of the crucial results in the alternative approach of Talagrand in [15] is a very general construction of pure states for measures in Hilbert spaces which provides a way around the representation results for exchangeable arrays that we are relying on here.

Since we will hereafter deal with the discrete directing measure (3.1), let us introduce some notation that will be more convenient throughout the rest of the paper. We will keep the dependence of the directing measure $\mu_\omega$ on $\omega$ implicit and write $\mathbb{E}$ for the integration in $\omega$. To describe an i.i.d. sample from measure $\mu$ in (3.1), consider i.i.d. random variables

$$\sigma^1, \sigma^2, \ldots \in \mathbb{N} \tag{3.3}$$

that take any value $l \geq 1$ with probability $w_l$, which is the weight corresponding to the index $l$ in the directing measure (3.1). Let us denote by $\langle \cdot \rangle$ the expectation in these random indices for a given measure $\mu$, that is, for any $n \geq 1$ and a function $h : \mathbb{N}^n \to \mathbb{R}$,

$$\langle h \rangle = \langle h(\sigma^1, \ldots, \sigma^n) \rangle = \sum_{l_1, \ldots, l_n} h(l_1, \ldots, l_n) w_{l_1} \cdots w_{l_n}. \tag{3.4}$$

By the configuration of $\mu$, we will understand the weights and configuration of its atoms,

$$w = (w_l)_{l \geq 1} \quad \text{and} \quad \mathcal{R} = (\xi(l) \cdot \xi(l'))_{l, l' \geq 1}. \tag{3.5}$$

Since $\xi(\sigma^l)$ are i.i.d. from distribution $\mu$, Proposition 1 can be rephrased by saying that the nondiagonal elements of the matrix $R$ can be generated by first generating a random measure $\mu$ [or its configuration (3.5)], then sampling indices (3.3) and setting

$$R_{l,l'} = \xi(\sigma^l) \cdot \xi(\sigma^{l'}) = \mathcal{R}_{\sigma^l, \sigma^{l'}}. \tag{3.6}$$

The Ghirlanda–Guerra identities (1.2) can be rewritten as

$$\mathbb{E}\langle f_n \psi(R_{1,n+1}) \rangle = \frac{1}{n} \mathbb{E}\langle f_n \rangle \mathbb{E}\langle \psi(R_{1,2}) \rangle + \frac{1}{n} \sum_{l=2}^{n} \mathbb{E}\langle f_n \psi(R_{1,l}) \rangle \tag{3.7}$$

and (1.1) can be rewritten as

$$\mathbb{E}\langle I(R_{1,2} = q_l) \rangle = m_{l+1} - m_l \quad \text{for } 1 \leq l \leq k. \tag{3.8}$$

The fact that in (3.1), all $\|\xi(l)\|^2 = q_k$ implies that in (2.3), all $a_l = 1 - q_k$ and, therefore, we can safely omit the term $a_l \delta_{l,l'}$ in (2.3) and redefine the matrix $R$ by $R_{l,l'} = \xi(\sigma^l) \cdot \xi(\sigma^{l'})$ for all $l, l' \geq 1$ so that, from now on, the diagonal elements are equal to $q_k$. As mentioned in the Introduction, as in [3], the crucial step which will allow us to use the Dovbysh–Sudakov representation (2.3) in order to make the induction step in the proof of Theorem 1 is the following.

THEOREM 3. *The matrix $\mathcal{R}$ in (3.5) is weakly exchangeable conditionally on $w = (w_l)$.*

Of course, this means that $\mathcal{R}$ is also weakly exchangeable unconditionally, which is how it will be used in the proof of Theorem 1, but the proof of the stronger statement of conditional exchangeability is exactly the same. The proof of Theorem 3 will be based on a certain invariance property of the joint distribution of $w$ and $\mathcal{R}$—quasi-stationarity under free evolution in the terminology of [3]—which, in our setting, will follow from the Ghirlanda–Guerra identities. Consider i.i.d. Rademacher random variables $(\varepsilon_l)_{l\geq 1}$ independent of the measure $\mu$. Given $t \geq 0$, consider a new sequence of weights

$$(3.9) \qquad w_l^t = \frac{w_l e^{t\varepsilon_l}}{\sum_{p\geq 1} w_p e^{t\varepsilon_p}},$$

defined by a random change of density proportional to $e^{t\varepsilon_l}$. Of course, these weights are not necessarily decreasing anymore, so let us denote by $(w_l^\pi)$ the weights $(w_l^t)$ arranged in decreasing order and let $\pi : \mathbb{N} \to \mathbb{N}$ be the permutation keeping track of where each index came from, $w_l^\pi = w_{\pi(l)}^t$. Let us define by

$$(3.10) \qquad \mu^\pi = \sum_{l\geq 1} w_l^\pi \delta_{\xi(\pi(l))} \quad \text{and} \quad \mathcal{R}^\pi = (\xi(\pi(l)) \cdot \xi(\pi(l')))_{l,l'\geq 1}$$

the probability measure $\mu$ after the change of density proportional to $e^{t\varepsilon_l}$ and the matrix $\mathcal{R}$ rearranged according to the reordering of weights. Analogously to (3.4), let us denote by $\langle h_n \rangle_\pi$ the average

$$(3.11) \qquad \langle h \rangle_\pi = \sum_{l_1,\ldots,l_n} h(l_1,\ldots,l_n) w_{l_1}^\pi \cdots w_{l_n}^\pi$$

and let $\mathbb{E}$ now denote the expectation in the randomness of the measure $\mu$ and the Rademacher sequence $(\varepsilon_l)$. Theorem 3 is a consequence of the following invariance principle.

THEOREM 4. *For $t < 1/2$, we have $(w^\pi, \mathcal{R}^\pi) \stackrel{\mathcal{D}}{=} (w, \mathcal{R})$.*

This result can be expressed by saying that the directing measure $\mu$ is stochastically stable under the change of density (3.9) and is a nontrivial consequence of the Ghirlanda–Guerra identities; in fact, as mentioned in the Introduction, this is, in some sense, a reversal of the usual derivation of the Ghirlanda–Guerra identities. There is a very important technical reason why we use Rademacher instead of the more obvious Gaussian change of density (as in [3]) that will become clear from the proof. Of course, as will be shown in the proof of Theorem 3, Theorem 4 for Rademacher change of density (3.9) implies the same result for Gaussian change of density as well.

Let us now show how Theorem 3 can be used to prove ultrametricity by induction on $k$.

PROOF OF THEOREM 1. Given a Gram–de Finetti matrix $R$, we consider the configuration matrix $\mathcal{R}$ of its directing measure defined in (3.5), which is symmetric and nonnegative definite. By (3.1), $\mathcal{R}_{l,l} = q_k$ and by (1.1), with probability 1, nondiagonal elements

$$\mathcal{R}_{l,l'} \in \{q_1, \ldots, q_{k-1}\}. \tag{3.12}$$

By Theorem 3, $\mathcal{R}$ is weakly exchangeable and, thus, is a Gram–de Finetti matrix. Therefore, using Proposition 1, there exists a random probability measure $\eta'$ on $H \times [0, \infty)$ such that

$$(\mathcal{R}_{l,l'}) \stackrel{\mathcal{D}}{=} (y_l \cdot y_{l'} + b_l \delta_{l,l'}), \tag{3.13}$$

where $(y_l, b_l)$ is an i.i.d. sequence from the distribution $\eta'$. Let us now show that (3.12) implies that if $\mu'$ is the marginal of $\eta'$ on $H$, then $\mu' = \sum_{l \geq 1} w'_l \delta_{\xi'(l)}$ for some distinct sequence $\xi'(l) \in H$, $w'_1 \geq w'_2 \geq \cdots > 0$, and, moreover,

$$\xi'(l) \cdot \xi'(l'), \qquad \|\xi'(l)\|^2 \in \{q_1, \ldots, q_{k-1}\}. \tag{3.14}$$

Indeed, if a point $y$ belongs to the support of $\mu'$, in the sense that $\mu'(B_\varepsilon(y)) > 0$ for all $\varepsilon > 0$, then in the i.i.d. sequence $(y_l)$ from this distribution, there will be infinitely many elements from $B_\varepsilon(y)$. Since, by (3.12), the scalar product $y_l \cdot y_{l'}$ of these elements belongs to $\{q_1, \ldots, q_{k-1}\}$ with probability 1, letting $\varepsilon \to 0$ proves that $\|y\|^2 \in \{q_1, \ldots, q_{k-1}\}$. If $z$ is another point in the support of $\mu'$ such that $z \in B_\varepsilon(y)$, then $\|z\|^2 \in \{q_1, \ldots, q_{k-1}\}$ and, therefore, $\|z\|^2 = \|y\|^2$ if $\varepsilon$ is small enough. The same argument also proves that $y \cdot z \in \{q_1, \ldots, q_{k-1}\}$, which implies that $y = z$ if $\varepsilon$ is small enough. This proves that the measure $\mu'$ is discrete and (3.14) holds. [*Remark:* Note that (3.14) guarantees that $q_{k-1} \geq 0$ and by a forthcoming induction, one can similarly conclude that all $q_l \geq 0$, which means that we did not need to invoke Talagrand's positivity principle and assume that $q_1 \geq 0$.]

Let us now explain the induction step. If $r^- = \min(r, q_{k-1})$, then the truncated matrix $R^- = (R^-_{l,l'})$ is symmetric, weakly exchangeable and, obviously, automatically satisfies the Ghirlanda–Guerra identities (1.2). Also, since $R_{l,l'} = \mathcal{R}_{\sigma^l, \sigma^{l'}}$ and $R^-_{l,l'} = \mathcal{R}^-_{\sigma^l, \sigma^{l'}}$, $R^-$ is nonnegative definite whenever $\mathcal{R}^-$ is, and the fact that $\mathcal{R}^-$ is nonnegative definite can be seen as follows. By (3.13), (3.14) and since $\mathcal{R}_{l,l} = q_k$ by (3.1), we get

$$(\mathcal{R}^-_{l,l'}) \stackrel{\mathcal{D}}{=} (y_l \cdot y_{l'} + (q_{k-1} - \|y_l\|^2)\delta_{l,l'}). \tag{3.15}$$

Since, by (3.14), $\|y_l\|^2 \leq q_{k-1}$, the right-hand side is obviously nonnegative definite. This implies that $\mathcal{R}^-$ and, therefore, $R^-$ are nonnegative definite and we have proven that the truncation $R^-$ is again a Gram–de Finetti matrix that satisfies the Ghirlanda–Guerra identities. Finally, the elements of $R^-$ take $k-1$ values $\{q_1,\ldots,q_{k-1}\}$,

$$(3.16) \qquad \mathbb{P}(R_{1,2}^- = q_{k-1}) = \mathbb{P}(R_{1,2} \geq q_{k-1}) = 1 - m_{k-1}$$

and, for $l \leq k-2$,

$$(3.17) \qquad \mathbb{P}(R_{1,2}^- = q_l) = \mathbb{P}(R_{1,2} = q_l) = m_{l+1} - m_l.$$

By the induction assumption, the matrix $R^-$ is ultrametric and together with (3.2), this completes the proof of Theorem 1. □

It is known [5] that an ultrametric matrix $R$ that satisfies the Ghirlanda–Guerra identities must be generated by the directing measure $\mu$ defined via the so called Derrida–Ruelle probability cascades [4, 11].

**4. Invariance and exchangeability.** It remains to prove the main consequences of the Ghirlanda–Guerra identities for the configuration of the directing measure $\mu$: a form of the stochastic stability of Theorem 4 and its application to the exchangeability of Theorem 3. In the proof below, we will need a couple of well-known properties of the Poisson–Dirichlet distribution $PD(s)$ which we will now recall.

By Theorem 2(c), the sequence $(w_l)$ in (3.1) has the Poisson–Dirichlet distribution $PD(s)$ with $s = m_k$ and it is defined by $w_l = u_l / \sum_{p \geq 1} u_p$, where $(u_l)_{l \geq 1}$ is the decreasing enumeration of a Poisson point process on $(0, \infty)$ with intensity measure $x^{-1-s}\, dx$. Let us consider an i.i.d. sequence $(X_l, Y_l)$ on $(0, \infty) \times \mathbb{R}$, independent of $(u_l)$ and such that $\mathbb{E} X_1 < \infty$. Note that $X_l$ and $Y_l$ need not be independent and, for example, $X_l$ can be a function of $Y_l$. Let $(Y_l')$ be an i.i.d. sequence independent of everything else such that for any measurable bounded function $\phi$,

$$(4.1) \qquad \mathbb{E}\phi(Y_1') = \frac{\mathbb{E} X_1^s \phi(Y_1)}{\mathbb{E} X_1^s},$$

which means that the distribution of $Y_1'$ is the distribution of $Y_1$ under the change of density $X_1^s / \mathbb{E} X_1^s$. Let $\theta$ be a random permutation of integers such that $(u_{\theta(l)} X_{\theta(l)})$ is arranged in decreasing order. The first property that will be useful (see Proposition 2.3 in [11] or Proposition 6.5.15 in [14]) states that

$$(4.2) \qquad (u_{\theta(l)} X_{\theta(l)}) \stackrel{\mathcal{D}}{=} ((\mathbb{E} X_1^s)^{1/s} u_l)$$

and, in particular, the sequence of weights $w_l' = u_l X_l / \sum_{p \geq 1} u_p X_p$, after rearranging in decreasing order, $(w_{\theta(l)}')$, again has the Poisson–Dirichlet distribution $PD(s)$. A more subtle result, Proposition A.2 in [4] (this result was rediscovered a couple of times—see Proposition 3.1 in [12] or Lemma 1.1 in [10]) implies that

$$(4.3) \qquad (u_{\theta(l)} X_{\theta(l)}, Y_{\theta(l)}) \stackrel{\mathcal{D}}{=} ((\mathbb{E} X_1^s)^{1/s} u_l, Y_l')$$

and, in particular, $(w'_{\theta(l)}, Y_{\theta(l)}) \stackrel{\mathcal{D}}{=} (w_l, Y_l')$. This property holds in more generality, but the case of real-valued $(Y_l)$ will be sufficient for our purposes.

PROOF OF THEOREM 4. Given $n \geq 2$, let us consider a function $f_n = f((R_{l,l'})_{1 \leq l < l' \leq n})$, where $R_{l,l'} = \xi(\sigma^l) \cdot \xi(\sigma^{l'})$, and suppose that $\|f\|_\infty \leq 1$. Consider a function $\varphi(t) = \mathbb{E}\langle f_n \rangle_\pi$. The central idea of the proof is to show that

$$(4.4) \qquad \varphi(t) = \varphi(0) \qquad \text{for } t < \frac{1}{2},$$

which means that a weakly exchangeable matrix $(R_{l,l'})_{l,l' \geq 1}$ has the same distribution under the directing measures $\mu^\pi$ and $\mu$. After we prove (4.4), we will explain how this implies that the configurations of random measures $\mu$ and $\mu^\pi$ have the same distribution, which is precisely the statement of the theorem. If, for $l \geq 2$, we let

$$(4.5) \qquad \Delta_l = \varepsilon(\sigma^1) + \cdots + \varepsilon(\sigma^{l-1}) - (l-1)\varepsilon(\sigma^l)$$

[we will write $\varepsilon(l)$ instead of $\varepsilon_l$], then it is easy to see that $\varphi'(t) = \mathbb{E}\langle f_n \Delta_{n+1} \rangle_\pi$ and, more generally,

$$\varphi^{(k)}(t) = \mathbb{E}\langle f_n \Delta_{n+1} \cdots \Delta_{n+k} \rangle_\pi.$$

Since $|\Delta_l| \leq 2l$, we get $|\varphi^{(k)}(t)| \leq 2^k (n+k)!/n!$. If we can show that $\varphi^{(l)}(0) = 0$ for all $l \geq 1$, then

$$|\varphi(t) - \varphi(0)| \leq \frac{(n+k)!}{k!n!} 2^k t^k$$

and letting $k \to +\infty$ implies that $\varphi(t) = \varphi(0)$ for $t < 1/2$. This is a good time to mention why we used a Rademacher sequence in the change of density (3.9) instead of the more obvious Gaussian. In the latter case, using Gaussian integration by parts, it is very easy to show that the Ghirlanda–Guerra identities imply $\varphi^{(l)}(0) = 0$; moreover, the Aizenman–Contucci identities [1] suffice here. However, the problem of controlling the derivatives $|\varphi^{(k)}(t)|$ becomes extremely difficult. Using bounded Rademacher random variables in the change of density gives us control of these derivatives for free, but it transfers the difficulty to showing that $\varphi^{(l)}(0) = 0$, which we now address. To complete the proof of (4.4), it remains to show that for all $k \geq 1$,

$$(4.6) \qquad \varphi^{(k)}(0) = \mathbb{E}\langle f_n \Delta_{n+1} \cdots \Delta_{n+k} \rangle = 0.$$

Since for $t = 0$, $\mu^\pi = \mu$ and, thus, $\langle \cdot \rangle$ is independent of the Rademacher sequence $(\varepsilon(l))$,

$$\varphi^{(k)}(0) = \mathbb{E}\langle f_n \mathbb{E}_\varepsilon \Delta_{n+1} \cdots \Delta_{n+k} \rangle,$$

where $\mathbb{E}_\varepsilon$ denotes the expectation in Rademacher random variables only. If $k$ is odd, then the derivative is zero by changing $(\varepsilon_l) \to (-\varepsilon_l)$. From now on, we will assume that $k$ is even. It is obvious that $\mathbb{E}_\varepsilon \Delta_{n+1} \cdots \Delta_{n+k}$ is the function of $(I(\sigma^l = \sigma^{l'}))_{l<l'}$ only and, by (3.1), $\sigma^l = \sigma^{l'}$ if and only if $R_{l,l'} = q_k$, which suggests that the Ghirlanda–Guerra identities can be used in the computation of (4.6). Let us start by expanding the product $\Delta_{n+1} \cdots \Delta_{n+k}$. Each term in the expansion corresponds to a collection $I$ of $n+k$ disjoint sets, $I_1, \ldots, I_{n+k}$, such that

$$\{n+1, \ldots, n+k\} = I_1 \cup \cdots \cup I_{n+k}$$

and such that $I$ describes the fact that we select each $\varepsilon(\sigma^j)$ for $1 \leq j \leq n+k$ from factors $\Delta_l$ with indices $l \in I_j$. We will call such a collection $I$ a *partition* of $\{n+1, \ldots, n+k\}$, even though some of the sets $I_1, \ldots, I_{n+k}$ can be empty. Therefore, we can write

(4.7) $$\Delta_{n+1} \cdots \Delta_{n+k} = \sum_I c_I \varepsilon(\sigma^1)^{|I_1|} \cdots \varepsilon(\sigma^{n+k})^{|I_{n+k}|}$$

for some constants $c_I$ that, of course, depend on the partitions $I$. Next, for any partition $P$ of $\{1, \ldots, n+k\}$, let us write

$$l \sim_P l' \iff l \text{ and } l' \text{ belong to the same element of } P$$

and let us denote by $I_P = I_P(\sigma^1, \ldots, \sigma^{n+k})$ the indicator of the event

(4.8) $$I_P = I\{\sigma^l = \sigma^{l'} \text{ if and only if } l \sim_P l', 1 \leq l, l' \leq n+k\}.$$

Using the fact that $1 = \sum_P I_P$, let us write, for any partition $I$ in (4.7),

$$\mathbb{E}_\varepsilon \varepsilon(\sigma^1)^{|I_1|} \cdots \varepsilon(\sigma^{n+k})^{|I_{n+k}|} = \sum_P I_P \mathbb{E}_\varepsilon \varepsilon(\sigma^1)^{|I_1|} \cdots \varepsilon(\sigma^{n+k})^{|I_{n+k}|}.$$

Each term on the right-hand side is either $I_P$ when

(4.9) $$\varepsilon(\sigma^1)^{|I_1|} \cdots \varepsilon(\sigma^{n+k})^{|I_{n+k}|} \equiv 1,$$

that is, when, for each set of the partition $P$, the number of factors $\varepsilon(\sigma^j)$ (with their multiplicities) with indices $j$ inside this set is even, or 0 otherwise, since, in this case, at least one independent factor $\varepsilon(\sigma^j)$ will remain and

$$\mathbb{E}_\varepsilon \varepsilon(\sigma^1)^{|I_1|} \cdots \varepsilon(\sigma^{n+k})^{|I_{n+k}|} = 0.$$

Therefore, if we denote by $\mathcal{P}(I)$ the collection of partitions $P$ of the first type for which (4.9) holds, then

$$(4.10) \qquad \mathbb{E}\langle f_n \mathbb{E}_\varepsilon \Delta_{n+1} \cdots \Delta_{n+k}\rangle = \sum_I c_I \sum_{P \in \mathcal{P}(I)} \mathbb{E}\langle f_n I_P\rangle.$$

Since $f_n$ is a function of the overlaps $(R_{l,l'})$ which takes only a finite number of values (1.1), it can be written as a linear combination of indicator functions of sets of the type

$$(4.11) \qquad \{R_{l,l'} = q_{l,l'} : 1 \leq l, l' \leq n\}$$

for any symmetric nonnegative definite matrix $(q_{l,l'})_{1 \leq l, l' \leq n}$ with $q_{l,l'} \in \{q_1, \ldots, q_k\}$ and diagonal elements $q_{l,l} = q_k$. Therefore, we can assume that $f_n$ is the indicator of the set (4.11). By (3.1) or by (3.2), we can assume that the constraints $(q_{l,l'})$ in (4.11) induce a partition $Q$ on the set $\{1, \ldots, n\}$ according to the rule

$$(4.12) \qquad l \sim_Q l' \quad \text{if and only if} \quad q_{l,l'} = q_k$$

and constraints $q_{l,l'}$ depend on $l, l'$ only through the partition elements in $Q$ which they belong to. If partition $Q$ consists of sets $Q_1, \ldots, Q_p$ and $l_j = \min\{l : l \in Q_j\}$ is the smallest index in each set, then $f_n$ can be written as

$$(4.13) \quad f_n = f'_n I_Q, \qquad \text{where } f'_n = I(\{R_{l,l'} = q_{l,l'} : l, l' \in \{l_1, \ldots, l_p\}\}).$$

In this representation, we separate constraints which describe how coordinates group together in the partition $Q$ from constraints between representatives $l_1, \ldots, l_r$ of each element of the partition, defined by $f'_n$. Note that in the definition of $f'_n$, all $q_{l,l'} \neq q_k$ for $l \neq l'$.

Returning to (4.10), if a partition $P$ of $\{1, \ldots, n+k\}$ does not agree with $Q$ on $\{1, \ldots, n\}$, then $f_n I_P = f'_n I_Q I_P \equiv 0$. This means that in (4.10), we can redefine $\mathcal{P}(I)$ to include only partitions $P$ that agree with $Q$, that is, $I_Q I_P = I_P$. For such partitions, we will now compute

$$\mathbb{E}\langle f_n I_P\rangle = \mathbb{E}\langle f'_n I_P\rangle$$

slightly more explicitly. Suppose that

$$P = P_1 \cup \cdots \cup P_p \cup P_{p+1} \cup \cdots \cup P_r$$

(we will abuse notation and write a partition as a union of its elements), where $P_l \cap \{1, \ldots, n\} = Q_l$ for $1 \leq l \leq p$ and $P_l \subseteq \{n+1, \ldots, n+k\}$ for $p < l \leq r$. Of course, it is possible that $r = p$. Our immediate goal will be to demonstrate that the Ghirlanda–Guerra identities imply that

$$(4.14) \qquad \mathbb{E}\langle f'_n I_P\rangle = \Phi(P)\mathbb{E}\langle f'_n\rangle$$

for some function $\Phi(P)$ that depends only on the configuration of the partition $P$. The exact formula for $\Phi(P)$ will not be important, but what *will* be important is to observe that it does not depend on the constraints in (4.13) that define $f'_n$. For $1 \leq j \leq r$, we denote by

(4.15) $$l_j = \min\{l : l \in P_j\}$$

the smallest index in the set $P_j$. Obviously, this definition agrees with the previous definition of $l_j$ for $Q_j$. Let us consider one of the sets in the partition that contains at least two points, for example, $P_r$. Let $l$ be the largest index in $P_r$ and let $P'$ be the restriction of the partition $P$ to the set $\{1,\ldots,n+k\} \setminus \{l\}$. We can then write $I_P = I_{P'}I(\sigma^{l_r} = \sigma^l)$. By (3.1), $\{\sigma^l = \sigma^{l'}\} = \{R_{l,l'} = q_k\}$ and we can treat $I(\sigma^l = \sigma^{l'})$ as a function of $R_{l,l'}$ when using the Ghirlanda–Guerra identities. Therefore, (3.7) implies that

$$\mathbb{E}\langle f'_n I_P\rangle = \frac{1-m_k}{n+k-1}\mathbb{E}\langle f'_n I_{P'}\rangle + \frac{1}{n+k-1}\sum_{j\neq l_r,l}\mathbb{E}\langle f'_n I_{P'}I(\sigma^{l_r} = \sigma^j)\rangle.$$

The only nonzero terms in the last sum correspond to $j \in P_r \setminus \{l_r, l\}$ and for such $j$, the constraint $\sigma^{l_r} = \sigma^j$ is already included in $P'$, so $I_{P'}I(\sigma^{l_r} = \sigma^j) = I_{P'}$ and we get

$$\mathbb{E}\langle f'_n I_P\rangle = \frac{1-m_k}{n+k-1}\mathbb{E}\langle f'_n I_{P'}\rangle + \frac{|P_r|-2}{n+k-1}\mathbb{E}\langle f'_n I_{P'}\rangle = \frac{|P_r|-1-m_k}{n+k-1}\mathbb{E}\langle f'_n I_{P'}\rangle.$$

Recursively, we can sequentially remove all coordinates with indices in $P_r$, except $\sigma^{l_r}$. If we consider the partition

$$P' = P_1 \cup \cdots \cup P_{r-1} \cup \{l_r\},$$

then

$$\mathbb{E}\langle f'_n I_P\rangle = \frac{(|P_r|-1-m_k)\cdots(|P_r|-(|P_r|-1)-m_k)}{(n+k-1)\cdots(n+k-(|P_r|-1))}\mathbb{E}\langle f'_n I_{P'}\rangle.$$

We can carry out the same computation on each of the partitions $P_1, \ldots, P_{r-1}$. As a result, if we consider the partition

(4.16) $$P' = \{l_1\} \cup \cdots \cup \{l_p\} \cup \{l_{p+1}\} \cup \cdots \cup \{l_r\}$$

and denote

$$\kappa_j = (|P_j|-1-m_k)\cdots(|P_j|-(|P_j|-1)-m_k)$$

for $1 \leq j \leq r$ and $\kappa_j = 1$ if $|P_j| = 1$, then

(4.17) $$\mathbb{E}\langle f'_n I_P\rangle = \frac{\kappa_r \cdots \kappa_1}{(n+k-1)\cdots(n+k-(|P_r|-1)-\cdots-(|P_1|-1))} \times \mathbb{E}\langle f'_n I_{P'}\rangle.$$

Finally, let us simplify $\mathbb{E}\langle f'_n I_{P'}\rangle$. If $p = r$, then $f'_n I_{P'} = f'_n$. If $p < r$, then we continue and consider a partition

(4.18) $$P'' = \{l_1\} \cup \cdots \cup \{l_p\} \cup \{l_{p+1}\} \cup \cdots \cup \{l_{r-1}\}.$$

Then $I_{P'} = I_{P''} - \sum_{j=1}^{r-1} I_{P''} I(\sigma^{l_r} = \sigma^{l_j})$ and, using the Ghirlanda–Guerra identities,

$$\mathbb{E}\langle f'_n I_{P'}\rangle = \mathbb{E}\langle f'_n I_{P''}\rangle - \sum_{j=1}^{r-1} \frac{1}{r-1} \mathbb{E}\langle f'_n I_{P''}\rangle(1-m_k) = m_k \mathbb{E}\langle f'_n I_{P''}\rangle.$$

Recursively, we can remove all coordinates $\sigma^{l_{p+1}}, \ldots, \sigma^{l_r}$ to get

(4.19) $$\mathbb{E}\langle f'_n I_{P'}\rangle = m_k^{r-p} \mathbb{E}\langle f'_n\rangle.$$

In this last term, we do not need to write the indicator of the partition $\{l_1\} \cup \cdots \cup \{l_p\}$ since these constraints are already contained in the definition of $f'_n$. Therefore, we have proven (4.14) with

$$\Phi(P) = \frac{\kappa_r \cdots \kappa_1}{(n+k-1)\cdots(n+k-(|P_r|-1)-\cdots-(|P_1|-1))} m_k^{r-p}$$

and equation (4.10) becomes

(4.20) $$\mathbb{E}\langle f_n \mathbb{E}_\varepsilon \Delta_{n+1} \cdots \Delta_{n+k}\rangle = \left(\sum_I c_I \sum_{P \in \mathcal{P}(I)} \Phi(P)\right) \mathbb{E}\langle f'_n\rangle.$$

It seems difficult to show algebraically that $\sum_I c_I \sum_{P \in \mathcal{P}(I)} \Phi(P) = 0$. However, as mentioned above, one can note that the computation leading to (4.20) depends on $f_n$ only through $I_Q$ in (4.13) since we only used the fact that partitions $P \in \mathcal{P}(I)$ should agree with $Q$ on $\{1, \ldots, n\}$. Therefore, (4.20) takes exactly the same form for $f_n = I_Q$, for which $f'_n$ is the indicator corresponding to the partition $Q_0 = \{l_1\} \cup \cdots \cup \{l_p\}$, that is,

$$\mathbb{E}\langle I_Q \mathbb{E}_\varepsilon \Delta_{n+1} \cdots \Delta_{n+k}\rangle = \left(\sum_I c_I \sum_{P \in \mathcal{P}(I)} \Phi(P)\right) \mathbb{E}\langle I_{Q_0}\rangle.$$

Another way to see this is simply to add up (4.20) for all $f_n$ corresponding to the same partition $Q$. Therefore, since $\mathbb{E}\langle I_{Q_0}\rangle \neq 0$, we will complete the proof if we can show that

$$\mathbb{E}\langle I_Q \mathbb{E}_\varepsilon \Delta_{n+1} \cdots \Delta_{n+k}\rangle = \mathbb{E}\langle I_Q \Delta_{n+1} \cdots \Delta_{n+k}\rangle = 0.$$

However, this is the $k$th derivative of the function $\varphi_Q(t) = \mathbb{E}\langle I_Q\rangle_\pi$ at $t = 0$ and the result will follow if we can show that $\varphi_Q(t) \equiv \varphi_Q(0)$. The crucial observation here is that $\mathbb{E}\langle I_Q\rangle_\pi$ depends only on the distribution of the sequence $(w_l^\pi)$ because

(4.21) $$\sigma^l = \sigma^{l'} \iff w_{\sigma^l}^\pi = w_{\sigma^{l'}}^\pi,$$

provided that all of the weights in $(w_l^\pi)$ are different with probability 1. By (4.2), $(w_l^\pi)$ has the Poisson–Dirichlet distribution $PD(m_k)$ and, therefore, all of the weights $w_l^\pi$ are different with probability 1, $(w_l^\pi)$ and $(w_l)$ have the same distribution and, by (4.21), $\mathbb{E}\langle I_Q\rangle_\pi = \mathbb{E}\langle I_Q\rangle$. This proves that $\varphi_Q(t) \equiv \varphi_Q(0)$ and completes the proof of (4.4).

If $R$ and $R^\pi$ are Gram–de Finetti matrices generated by the random directing measures $\mu$ and $\mu^\pi$, respectively, then (4.4) obviously implies that $R \stackrel{\mathcal{D}}{=} R^\pi$. It remains to prove that

$$R \stackrel{\mathcal{D}}{=} R^\pi \quad \Longrightarrow \quad (w, \mathcal{R}) \stackrel{\mathcal{D}}{=} (w^\pi, \mathcal{R}^\pi).$$

This will follow from the fact that, conditionally on $(w, \mathcal{R})$, the matrix $R$ is generated by $R_{l,l'} = \mathcal{R}_{\sigma^l, \sigma^{l'}}$, as in (3.6), from which one can show that $(w, \mathcal{R}) = \phi(R)$ almost surely for some measurable function $\phi$, that is, the configuration $(w, \mathcal{R})$ of the directing measure can be uniquely reconstructed from the overlap matrix $R$. Note that, by (3.1), with probability 1, the matrix $R$ is *ultrametric at the level* $k$, in the sense that the relation $l \sim_k l'$ defined by $R_{l,l'} = q_k$ is an equivalence relation on $\mathbb{N}$, and for any two equivalence classes $N_1$ and $N_2$, the coordinates $R_{l,l'}$ are equal for all $l \in N_1$ and $l' \in N_2$. Let $w^n(R)$ be the vector of frequencies of the equivalence classes restricted to the set $\{1, \ldots, n\}$, arranged in decreasing order and then extended to an infinite vector by appending all zeros. Let $\mathcal{R}^n(R)$ be the matrix of overlaps between the equivalence classes defined by

$$(\mathcal{R}^n(R))_{l,l'} = R_{i,j}$$

for any representatives $i$ and $j$ of the equivalence classes corresponding to nonzero $w^n(l)$ and $w^n(l')$ and extended to an infinite matrix by setting $q_k$ on the diagonal and zeros everywhere else. Define $\phi(R)$ as the coordinate-wise limit

$$\phi(R) = \lim_{n \to +\infty} (w^n(R), \mathcal{R}^n(R)) \tag{4.22}$$

if such limit exists and some fixed value otherwise, including when the matrix $R$ is not ultrametric at the level $k$. Since, conditionally on $(w, \mathcal{R})$, the matrix $R$ is generated as in (3.6), by the strong law of large numbers, $w^n(R)$ converges almost surely to $w$ coordinate-wise. All coordinates of $w$ are different with probability 1 since the sequence $(w_l)$ has the Poisson–Dirichlet distribution. For any such $w$, given $(w, \mathcal{R})$, $\mathcal{R}^n(R)$ also converges almost surely to $\mathcal{R}$ coordinate-wise because, asymptotically, each equivalence class corresponds to a unique $\xi(i)$ in the support of $\mu$. This proves that $(w, \mathcal{R}) = \phi(R)$ with probability 1 and since, by (4.4), the distribution of $R$ is the same under the directing measures $\mu$ or $\mu^\pi$, we get $(w^\pi, \mathcal{R}^\pi) \stackrel{\mathcal{D}}{=} (w, \mathcal{R})$. □

Finally, it remains to prove that the invariance principle of Theorem 4 implies exchangeability of the matrix $\mathcal{R}$.

PROOF OF THEOREM 3. Even though the underlying idea of our proof is the same as the idea in Proposition 3.3 in [3], we provide a more direct argument based on a very explicit control of the mixing induced by the random change of density (3.9). Let us start with the observation that Theorem 4 also holds if we replace a Rademacher sequence in the change of density (3.9) by an i.i.d. standard Gaussian sequence $(g_l)$.

LEMMA 2. *Theorem 4 holds with the change of density*

$$w_l^t = \frac{w_l e^{tg_l}}{\sum_{p \geq 1} w_p e^{tg_p}} \tag{4.23}$$

*for arbitrary $t > 0$ and i.i.d. standard Gaussian $(g_l)$.*

PROOF. This follows from the fact that the invariance provided by Theorem 4 will be preserved if we iterate the process of making the change of density (3.9). Namely, if $(\varepsilon_l^k)_{l \geq 1}$ are i.i.d. copies of $(\varepsilon_l)_{l \geq 1}$ for $k \geq 1$ and if we define $S_l^k = \varepsilon_l^1 + \cdots + \varepsilon_l^k$, then replacing (3.9) with

$$v_l^t = \frac{w_l e^{tS_l^k}}{\sum_{p \geq 1} w_p e^{tS_p^k}}, \tag{4.24}$$

the statement of Theorem 4 still holds. We can replace $t$ in (4.24) by $tk^{-1/2}$ for any fixed $t > 0$ and large enough $k$ such that $tk^{-1/2} < 1/2$. Each element of the i.i.d. sequence $(k^{-1/2}S_l^k)_{l \geq 1}$ converges in distribution to the standard Gaussian as $k \to +\infty$. Therefore, we can choose these sequences for all $k \geq 1$ on the same probability space with some i.i.d. Gaussian sequence $(g_l)_{l \geq 1}$ so that $k^{-1/2}S_l^k \to g_l$ almost surely for each $l$. It is easy to check that the sum

$$\sum_{p \geq 1} w_p \exp(tk^{-1/2}S_p^k) \to \sum_{p \geq 1} w_p \exp tg_p \qquad \text{a.s.,}$$

possibly over some subsequence $(k(n))$. Then $(v_l^t)$ converges almost surely to the sequence (4.23) and, as a result, $(w^\pi, \mathcal{R}^\pi)$ defined in terms of (4.24) also converges almost surely to the corresponding configuration defined in terms of $(g_l)$. Since, by Theorem 4, the distribution of $(w^\pi, \mathcal{R}^\pi)$ remains the same along this sequence, the distribution of the limiting configuration $(w^\pi, \mathcal{R}^\pi)$ defined in terms of $(g_l)$ is equal to the distribution of $(w, \mathcal{R})$. □

Now, given $n \geq 1$, let us consider a fixed permutation $\rho$ of $\{1, \ldots, n\}$ and a measurable subset $A \subseteq [-1, 1]^{n^2}$. Given $m \geq 1$, consider a measurable subset $B \subseteq [0, 1]^m$. To prove that conditionally on $w$ the overlap matrix $\mathcal{R}$ is weakly exchangeable, we need to show that

$$\mathbb{P}(\mathcal{R}_n^\rho \in A, (w_l)_{l \leq m} \in B), \qquad \text{where } \mathcal{R}_n^\rho = (\mathcal{R}_{\rho(l), \rho(l')})_{1 \leq l, l' \leq n},$$

does not depend on the permutation $\rho$. Without loss of generality, we can assume that $m = n$, by redefining the sets $A$ and $B$ and permutation $\rho$. Also, to simplify notation, we will write $w \in B$ instead of $(w_l)_{l \leq m} \in B$. Let $\pi$ be a permutation of indices induced by the rearrangement of the sequence (4.23), that is, $w_l^\pi = w_{\pi(l)}^t$. If we let

$$\mathcal{R}_n^{\pi \circ \rho} = (\mathcal{R}_{\pi \circ \rho(l), \pi \circ \rho(l')})_{1 \leq l, l' \leq n},$$

then Theorem 4 implies that

(4.25) $\quad \mathbb{P}(\mathcal{R}_n^\rho \in A, w \in B) = \mathbb{P}(\mathcal{R}_n^{\pi \circ \rho} \in A, w^\pi \in B).$

Intuitively, when $t$ goes to infinity, the order of $\pi(1), \ldots, \pi(n)$ becomes completely random because it is determined by the order of $\log w_{\pi(l)} + t g_{\pi(l)}$ for $1 \leq l \leq n$, which is asymptotically, for $t \to +\infty$, determined by the order of $g_{\pi(1)}, \ldots, g_{\pi(n)}$. Therefore, in the limit, the distribution of $\mathcal{R}_n^{\pi \circ \rho}$, and, thus, of $\mathcal{R}_n$, should not depend on $\rho$, which means that $\mathcal{R}$ is weakly exchangeable. However, since, a priori, we do not control the dependence of $w$ and $\mathcal{R}$, turning this intuition into a rigorous argument requires some work. Let us denote by $j = (j(1), \ldots, j(n))$ a generic vector with all indices $j(l)$ different and let $\pi \circ \rho = (\pi \circ \rho(1), \ldots, \pi \circ \rho(n))$. With this notation, the right-hand side of (4.25) can be written as

(4.26)
$$\sum_j \mathbb{P}(\mathcal{R}_n^{\pi \circ \rho} \in A, w^\pi \in B, \pi \circ \rho = j)$$
$$= \sum_j \mathbb{P}(\mathcal{R}_n^j \in A, w^\pi \in B, \pi \circ \rho = j),$$

where we have also introduced the notation $\mathcal{R}_n^j = (\mathcal{R}_{j(l), j(l')})_{1 \leq l, l' \leq n}$. Conditionally on $w = (w_l)$, the events $\{\mathcal{R}_n^j \in A\}$ and $\{w^\pi \in B, \pi \circ \rho = j\}$ are independent since the latter depends only on the sequence $(g_l)$ and, therefore,

(4.27)
$$\mathbb{P}(\mathcal{R}_n^j \in A, w^\pi \in B, \pi \circ \rho = j | w)$$
$$= \mathbb{P}(\mathcal{R}_n^j \in A | w) \mathbb{P}(w^\pi \in B, \pi \circ \rho = j | w).$$

If $\tau$ is another fixed permutation of $\{1, \ldots, n\}$, then (4.25), (4.26) and (4.27) imply that

(4.28)
$$|\mathbb{P}(\mathcal{R}_n^\rho \in A, w \in B) - \mathbb{P}(\mathcal{R}_n^\tau \in A, w \in B)|$$
$$\leq \sum_j \int |\mathbb{P}(w^\pi \in B, \pi \circ \rho = j | w)$$
$$- \mathbb{P}(w^\pi \in B, \pi \circ \tau = j | w)| \, d\Lambda(w),$$

where $\Lambda$ is the distribution of $w$. Let us express one of the events $C_\rho = \{w^\pi \in B, \pi \circ \rho = j\}$ in terms of the sequence $(g_l)$. If we let

(4.29) $\quad k = j \circ \rho^{-1}, \quad$ that is, $\quad k(l) = j(\rho^{-1}(l)) \qquad$ for $1 \leq l \leq n,$

then, by the definition of $\pi$, the event $\{\pi \circ \rho = j\}$ expresses the fact that for $1 \leq l \leq n$, the number $w_{k(l)} \exp t g_{k(l)}$ occupies the position $l$ among all the elements of $(w_i \exp t g_i)$ arranged in decreasing order. If we introduce the notation

$$\gamma_{k(l)} = t^{-1} \log w_{k(l)}, \qquad z_l = g_{k(l)} + \gamma_{k(l)},$$
$$x = \sup_{i \notin j}(g_i + t^{-1} \log w_i), \qquad y = \sum_{i \notin j} w_i e^{tg_i},$$

then the event $C_\rho$ can be written as

(4.30) $\quad C_\rho = \left\{\left(\frac{e^{tz_l}}{y + \sum_{1 \leq i \leq n} e^{tz_i}}\right)_{1 \leq l \leq n} \in B, z_1 \geq \cdots \geq z_n \geq x\right\}.$

Let us first consider the probability of $C_\rho$ conditionally on $w$ and $(g_i)_{i \notin j}$, that is, for fixed $x, y$ and $(\gamma_{k(l)})$. Since $(z_l)$ are independent and $z_l$ has normal distribution $N(\gamma_{k(l)}, 1)$, we can write

$$\mathbb{P}(C_\rho | w, (g_i)_{i \notin j}) = \frac{1}{(\sqrt{2\pi})^n} \int_{C_\rho} \exp\left(-\frac{1}{2} \sum_{l=1}^n (z_l - \gamma_{k(l)})^2\right) dz_1 \cdots dz_n$$

$$= \frac{1}{(\sqrt{2\pi})^n} \exp\left(-\frac{1}{2} \sum_{l=1}^n \gamma_{j(l)}^2\right)$$

$$\times \int_{C_\rho} \exp\left(\sum_{l=1}^n \gamma_{k(l)} z_l - \frac{1}{2} \sum_{l=1}^n z_l^2\right) dz_1 \cdots dz_n.$$

Since the event (4.30) does not explicitly depend on $\rho$, the last integral depends on $\rho$ only through the term $\sum_{l=1}^n \gamma_{k(l)} z_l$. If we denote by $(\gamma_l^-)$ and $(\gamma_l^+)$ the sequence $(\gamma_{k(l)})$ arranged in decreasing and increasing order, respectively, then, on the event $C_\rho$ (since $z_1 \geq \cdots \geq z_n$),

$$\sum_{l=1}^n \gamma_l^+ z_l \leq \sum_{l=1}^n \gamma_{k(l)} z_l \leq \sum_{l=1}^n \gamma_l^- z_l.$$

Therefore, the probability $\mathbb{P}(C_\rho | w, (g_i)_{i \notin j})$ is maximized on the permutation $\rho$ for which the sequence $k(l) = j(\rho^{-1}(l))$ in (4.29) is increasing, that is, $\rho$ and $j$ are similarly ordered. This is, obviously, equivalent to $\pi \circ e = j^+$, where $e$ is the identity permutation, $e(l) = l$, and $j^+$ is the increasing rearrangement of $j$. Similarly, $\mathbb{P}(C_\rho | w, (g_i)_{i \notin j})$ is minimized on the permutation $\rho$ for which the sequence in (4.29) is decreasing, which is equivalent to $\pi \circ e' = j^+$ for the inverse permutation, $e'(l) = n - l + 1$. Averaging over $(g_i)_{i \notin j}$, we have proven that

$$\mathbb{P}(w^\pi \in B, \pi \circ e' = j^+|w) \leq \mathbb{P}(w^\pi \in B, \pi \circ \rho = j|w)$$
$$\leq \mathbb{P}(w^\pi \in B, \pi \circ e = j^+|w)$$

and, therefore,

$$|\mathbb{P}(w^\pi \in B, \pi \circ \rho = j|w) - \mathbb{P}(w^\pi \in B, \pi \circ \tau = j|w)|$$
$$\leq \mathbb{P}(w^\pi \in B, \pi \circ e = j^+|w) - \mathbb{P}(w^\pi \in B, \pi \circ e' = j^+|w)$$

for any $\rho, \tau$ and $j$. Plugging this into (4.28) gives

(4.31)
$$\frac{1}{n!}|\mathbb{P}(\mathcal{R}_n^\rho \in A, w \in B) - \mathbb{P}(\mathcal{R}_n^\tau \in A, w \in B)|$$
$$\leq \mathbb{P}(w^\pi \in B, \exists j : \pi \circ e = j^+) - \mathbb{P}(w^\pi \in B, \exists j : \pi \circ e' = j^+).$$

We divide by $n!$ because each $j^+$ corresponds to $n!$ different $j$. It remains to show that the right-hand side goes to zero when $t$ in (4.23) goes to infinity. Let us recall the definition $w_l^\pi = w_{\pi(l)}^t$ and let us similarly define $g_l^\pi = g_{\pi(l)}$. The event $\{\exists j : \pi \circ e = j^+\}$ can then be expressed in terms of $(w^\pi, g^\pi)$, as follows. On one hand, this event simply means that $\pi(1) < \cdots < \pi(n)$. On the other hand, (4.23) implies that, if we let $\kappa = \sum_{p \geq 1} w_p e^{tg_p}$,

$$w_{\pi(l)} = \kappa w_{\pi(l)}^t e^{-tg_{\pi(l)}} = \kappa w_l^\pi e^{-tg_l^\pi}$$

and, therefore,

(4.32) $$\{\exists j : \pi \circ e = j^+\} = \{w_1^\pi e^{-tg_1^\pi} > \cdots > w_n^\pi e^{-tg_n^\pi}\}.$$

By (4.2), $(w_l^\pi)$ has the Poisson–Dirichlet distribution $\Lambda = PD(m_k)$ and it is easy to check that (4.3) and (4.1) imply that $(g_l^\pi)$ is an i.i.d. sequence with normal distribution $\nu = N(tm_k, 1)$, independent of $(w_l^\pi)$. Therefore,

$$\mathbb{P}(w^\pi \in B, \exists j : \pi \circ e = j^+)$$
$$= \int_B \nu^{\otimes n}((g_l^\pi)_{1 \leq l \leq n} : w_1^\pi e^{-tg_1^\pi} > \cdots > w_n^\pi e^{-tg_n^\pi}) d\Lambda(w^\pi).$$

For any fixed $w^\pi$, $\nu^{\otimes n}(w_1^\pi e^{-tg_1^\pi} > \cdots > w_n^\pi e^{-tg_n^\pi}) \to 1/n!$ when $t \to +\infty$ since, asymptotically, this event is equivalent to $g_1^\pi < \cdots < g_n^\pi$ and, therefore, $\mathbb{P}(w^\pi \in B, \exists j : \pi \circ e = j^+) \to \Lambda(B)$ as $t \to +\infty$. Similarly, the fact that $\mathbb{P}(w^\pi \in B, \exists j : \pi \circ e' = j^+) \to \Lambda(B)$, together with (4.31), completes the proof. □

**Acknowledgments.** The author would like to thank Michel Talagrand and the referees for many valuable comments and suggestions.

## REFERENCES

[1] AIZENMAN, M. and CONTUCCI, P. (1998). On the stability of the quenched state in mean-field spin-glass models. *J. Statist. Phys.* **92** 765–783. MR1657840

[2] ALDOUS, D. J. (1985). Exchangeability and related topics. In *École D'été de Probabilités de Saint-Flour, XIII—1983. Lecture Notes in Math.* **1117** 1–198. Springer, Berlin. MR883646

[3] ARGUIN, L.-P. and AIZENMAN, M. (2009). On the structure of quasi-stationary competing particles systems. *Ann. Probab.* **37** 1080–1113.

[4] BOLTHAUSEN, E. and SZNITMAN, A.-S. (1998). On Ruelle's probability cascades and an abstract cavity method. *Comm. Math. Phys.* **197** 247–276. MR1652734

[5] BOVIER, A. and KURKOVA, I. (2004). Derrida's generalised random energy models. I. Models with finitely many hierarchies. *Ann. Inst. H. Poincaré Probab. Statist.* **40** 439–480. MR2070334

[6] DOVBYSH, L. N. and SUDAKOV, V. N. (1982). Gram-de Finetti matrices. *Zap. Nauchn. Sem. Leningrad. Otdel. Mat. Inst. Steklov. (LOMI)* **119** 77–86, 238, 244–245. MR666087

[7] GHIRLANDA, S. and GUERRA, F. (1998). General properties of overlap probability distributions in disordered spin systems. Towards Parisi ultrametricity. *J. Phys. A* **31** 9149–9155. MR1662161

[8] PARISI, G. (1980). A sequence of approximate solutions to the S-K model for spin glasses. *J. Phys. A* **13** L115–L121.

[9] PANCHENKO, D. (2007). A note on Talagrand's positivity principle. *Electron. Comm. Probab.* **12** 401–410 (electronic). MR2350577

[10] PANCHENKO, D. and TALAGRAND, M. (2007). On one property of Derrida–Ruelle cascades. *C. R. Math. Acad. Sci. Paris* **345** 653–656. MR2371485

[11] RUELLE, D. (1987). A mathematical reformulation of Derrida's REM and GREM. *Comm. Math. Phys.* **108** 225–239. MR875300

[12] RUZMAIKINA, A. and AIZENMAN, M. (2005). Characterization of invariant measures at the leading edge for competing particle systems. *Ann. Probab.* **33** 82–113. MR2118860

[13] SHERRINGTON, D. and KIRKPATRICK, S. (1975). Solvable model of a spin glass. *Phys. Rev. Lett.* **35** 1792–1796.

[14] TALAGRAND, M. (2003). *Spin Glasses: A Challenge for Mathematicians. Ergebnisse der Mathematik und Ihrer Grenzgebiete. 3. Folge. A Series of Modern Surveys in Mathematics [Results in Mathematics and Related Areas. 3rd Series. A Series of Modern Surveys in Mathematics]* **46**. Springer, Berlin. MR1993891

[15] TALAGRAND, M. (2009). Construction of pure states in mean-field models for spin glasses. *Probab. Theory Related Fields*. To appear.

DEPARTMENT OF MATHEMATICS
TEXAS A&M UNIVERSITY
MAILSTOP 3386
COLLEGE STATION, TEXAS 77843
USA
E-MAIL: panchenk@math.tamu.edu